\documentclass[smallextended,numbook,runningheads]{svjour3}     % onecolumn (second format)
\smartqed  % flush right qed marks, e.g. at end of proof
\usepackage{graphicx}
\usepackage{amsmath}
\usepackage{mathptmx}
\usepackage{amssymb}
\usepackage{amsfonts,euscript}
\usepackage{amsmath}
\usepackage{color}
\journalname{BIT}

\newcommand {\RR}{\mathbb R}
\newcommand {\NN}{\mathbb N}

\newcommand {\II} {\mathbb I}
\newcommand {\CI} {\bar{\mathbb I}}
\newcommand {\OO}{\Omega}
\newcommand {\CO}{\overline\Omega}
\newcommand {\PO}{\partial\Omega}
\newcommand {\vx}{{\bf x}}
\newcommand {\vv}{{\bf v}}
\newcommand {\vd}{{\bf d}}

\begin{document}

\title{On the $L_p$-error of approximation of bivariate functions by harmonic splines}
%\vspace{2ex}

\author{Yuliya Babenko        \and
        Tatyana Leskevich %etc.
}

\institute{Y. Babenko \at
              Department of Mathematics and Statistics,
Kennesaw State University, Kennesaw, GA, USA, 30344\\
              \email{ybabenko@kennesaw.edu}           %  \\
           \and
           T. Leskevich \at
              Department of Mathematical Analysis and Theory of Functions,
Dnepropetrovsk National University,
Dnepropetrovsk, UKRAINE, 49050 \\
\email{tleskevich@gmail.com}
}

\date{Received: date / Accepted: date}

\maketitle

\begin{abstract}
Interpolation by various types of splines is the standard procedure
in many applications. In this paper we shall discuss harmonic spline
``interpolation'' (on the lines of a grid) as an alternative
to polynomial spline interpolation (at vertices of a grid). We will discuss some advantages
and drawbacks of this approach and present the asymptotics of the $L_p$-error for adaptive approximation by harmonic splines.
%Include keywords and mathematical subject classification numbers as needed.
\keywords{interpolation \and adaptive \and harmonic spline \and optimal error \and asymptotics}
\subclass{41A15 \and 41A60 \and 68W25 \and 97N50}
\end{abstract}

\section{Harmonic splines as an alternative to polynomial splines} \label{S1}

In ~\cite{nasha} authors investigated the question of adaptive
approximation by polynomial splines on box partitions in $\RR^d$.
The obtained a general description as well as sharp constant in
several popular cases. However, polynomial splines might not be the
best option when it comes to interpolating a function over box
partitions. There is an alternative - harmonic splines - which can
be viewed as a direct multivariate generalization of linear splines
in one dimension. This paper addresses approximation of bivariate
functions and, therefore, from now on we shall consider harmonic
splines of two variables.

The following several observations help to see the analogy between univariate
linear splines and bivariate harmonic splines on rectangular partitions more clear:

\begin{enumerate}
\item
Linear univariate splines interpolate the given function on the
whole boundary of the mesh which is the points of the partition of
the domain. While bilinear splines defined on
rectangular partitions interpolate the given function only at the
vertices of the mesh elements, not on the whole boundary, harmonic
splines are constructed to coincide with the original function on
the boundary of the meshes - on the lines (as oppose to just at the
points). Litvin in ~\cite{Litvin} introduces different terminology
to distinguish between these types of interpolation and he says that
harmonic splines ``inter{\it lineates}'' function instead of just
inter{\it polat}ing. Using harmonic splines has advantages and
drawbacks and the choice depends on a particular problem. The strong
advantage is that it uses more information about the function.
However, it also imposes restrictions on the mesh. In the case of
square domain (considered in this paper) the meshes will consist of
mostly squares. In the case of a rectangular domain, the elements of
the mesh will have the same shape (same width to length ratio) as
the domain, and only size will adapt to the local behavior of the
function.

\item Both univariate linear spline $P_1(x)$ and bivariate harmonic spline $U(x,y)$ satisfy corresponding operator equations. Indeed, if we denote by $\Delta_1=\frac{d^2}{dx^2}$ and $\Delta_2=\frac{\partial^2}{\partial x^2}+\frac{\partial^2}{\partial y^2}$, then clearly
$$
\Delta_1 P_1(x) =0\;\;\;\;\hbox{and}\;\;\;\; \Delta_2 U(x,y)=0.
$$
\item
A univariate linear interpolating spline and bivariate harmonic
spline attain maximum and minimum values on the boundary of the
element.
\item
If the values of the {univariate} linear spline at internal nodes
are determined from the continuity condition on the first
derivative, spline degenerates into a segment of straight line. If
the values of harmonic spline are determined from the condition of
continuity of the normal derivative to the interior boundary of the
mesh, then the spline degenerates into a harmonic function on the
union of domains (which is, by the way, far less restrictive).
\item Both have very similar local estimates of the error of approximation ~\cite{Klim}.
\end{enumerate}

For more discussion of properties of harmonic and polyharmonic
splines (or polysplines) see the monographs ~\cite{Koun,Litvin} and
references therein.

%\tc {mozhet tut dobavit' chto-to tipa: Esche odnim preimuschestvom
%garmonicheskix splaynov yavlyaetsya to, chto v otlichii ot ...}
%Unlike multivariate polynomials, polyharmonic functions,\tc{mozhet
%tut konkretnee o bilineynyx i garmonicheckix funkciyax} satisfy the
%generalized definition of a Chebyshev system and therefore create a
%basis for multivariate \tc{bivariate} interpolation theory.

To approximate the bivariate function (so that the approximant also
coincides with the given function on some lines) polynomial splines,
in particular, blending splines, can {also} be used. However, there
exists a broad class of functions that can be well approximated by
harmonic splines, but are not at all approximated by such polynomial
splines under the same requirements concerning the order of
smoothness. Let us present an example of such a function.

{\bf Example (~\cite{Klim}).} Let us consider the function
\begin{multline}
f(x,y)=\frac 4 {\pi} \displaystyle \sum_{k=0}^{\infty}\frac 1
{(2k+1)\alpha^{2}_k \sinh(\alpha_k b)}\left(
\sinh(\alpha_k(y-b))+\sinh(\alpha_k b)\right.\\ \left.
-\sinh(\alpha_k y) \right) \sin(\alpha_k x),\nonumber
\end{multline}
$$
0\leq x\leq a,\;0\leq y\leq b,\;\alpha_k=(2k+1)\pi/a.
$$
For this function
$$
\displaystyle \lim_{x\to 0, y\to 0} \frac{\partial ^2 f(x,y)}{\partial x \partial y} = \infty,
$$
and, hence, it cannot be approximated by blending polynomial
splines. On the other hand, ${\Delta f(x,y)=-1}$ which implies
that $f$ can be well approximated by harmonic splines.

In univariate case the idea of using splines constructed based on
solution of some differential operator equations, called
$L$-splines, has been more or less well-developed (see, for
example, ~\cite{Schum} and references therein). There has been
several attempts to generalize this concept to multivariate case. For example, Litivn in his monograph
~\cite{Litvin} introduces concepts of ``interlineation'' and
``interflatation''. These concepts are natural and direct
generalizations of interpolation to the case when one needs to
reconstruct the given function based on knowing values of function  (and possibly its normal
derivatives up to some order) on one or several lines
(``interlineation'') or $m$-dimensional $(1<m<d)$ linear manifolds (``interflatation'') in $\RR^d$.

As already mentioned earlier, in ~\cite{Litvin}  for the given on $\RR^2$ function Litvin suggests and compares
two methods of polynomial and non-polynomial ``interlineation'' on
several lines. The first (polynomial) is a generalization of Hermite
polynomials which preserves the
smoothness order of the original function. The second (not
polynomial) operator is a generalization of D'Alambert operator,
which is the solution of homogeneous wave equation, to the case when
data is given as values of the normal derivatives of all orders up to and including
 $N$ $(N>1)$.

%\tc{Ty prava, sleduyuschie dal'she primeneniya nuzhno sokrat', no me
%ochen' slozhno pravit' angliyskiy tekst, esli skazhesh, to ya
%luchshe po russki napishu kak ya eto vizhu, no bystro eto ne budet
%:)}

There exist many applications in which data is given not at a
discrete set of points but on some linear manifolds (lines, planes, etc.). One of
many examples is cartography. The measurement and analysis of
bathymetric measurements is one of the core areas of modern
hydrography, and a fundamental component in ensuring the safety of
water transportation, efficiency of offshore resource development,
precision and effectiveness in recovering underwater objects as well
as search and rescue missions. The data used to make bathymetric maps today typically
comes from an echosounder (sonar) mounted beneath or over the side
of a boat, ``pinging'' a beam of sound downward at the seafloor or
from remote sensing LIDAR or LADAR systems.
%The amount of time it takes for the sound or light to travel through the water, bounce off the seafloor, and return to the sounder tells the equipment what the distance to the seafloor is.
%Depending on the depth and complexity of the ocean floor, several types of sonars are used to survey and collect data, multibeam echosounder and side scan sonars are among them.
%Multibeam sonar measures the depth of the sea floor by analyzing the time it takes for sound waves to travel from a boat to the sea floor and back. Multibeam sonar is very useful in areas such as the Northeast U.S. and Alaska, where the sea floor is complex and often strewn with thousands of rocks.
In all the cases, the data obtained are univariate - along the lines, curves or, more generally, a grid composed of them - depending on the course of the surveyor vessel (which can be adapted when needed).
The goal is to recover the function describing the relief of the seafloor based on the ``traces'' of this function on the lines which described the course of the vessel with sonar.

%The detection, classification and localization performance of a sonar depends on the environment (water temperature, tides, etc.) and the receiving equipment, as well as the transmitting equipment in an active sonar or the target radiated noise in a passive sonar. Therefore, it is natural to assume that the data contains noise and measurement errors.

Similar technology is used for mapping surfaces of cosmic objects. Only instead of sonar, which uses sound waves, radars are used. Radars are based on high-frequency electromagnetic radiation which can penetrate the surroundings of objects, for instance, the clouds covering Venus.

Another example of applications to surface reconstruction based on data given on
lines, curves, or hyperplanes would be manufacturing of car,
aircraft etc. bodies. Currently, most popular tool used in these models is polynomial splines. Being easy to manipulate, they have many drawbacks: lack of global smoothness (if working with convenient low degree polynomial pieces) and lack of flexibility in choosing the mesh (for instance, it might be necessary to base a part of the body on hexagon or pentagon, not only triangle or rectangle). Harmonic and their generalization - polyharmonic - splines have great advantage in both directions: they preserve higher smoothness and they can be constructed on any (more or less reasonable) domain.

%Usually the body of an
%object is composed of many pieces built separately and then put
%together. Clearly, it is crucial that the pieces are put together
%smoothly and so that the aerodynamical qualities of the final
%product are maximized.
%
%
%
%

\section {Notation, statement of the main problem, and the main result} \label{S2}

Let in the space $\RR^2$ of points $\vx=(x_1,x_2)$ the unit square
$[0,1]\times[0,1]$, be given with interior denoted by $\II^2$, and
the closure denoted by $\CI^2$. By $\square_N$ we shall denote a
partition of $\CI^2$ whose number of elements has order $N$ as $N\to \infty$, and is so that the majority of elements are
squares, however some small number ($o(N)$ as $N\to \infty$) of rectangles is also allowed.
The interior of an arbitrary element of the partition $\square_N$ we shall denote by $\OO^N$, the closure by
$\CO^N$, and the boundary by $\PO^N$.

We shall need the following standard notation:
 $C(\OO)$ denotes the space of continuous in some region
$\OO$ functions, $C^2(\OO)$ denotes the space of twice
differentiable in $\OO$ functions, and $L_p(\OO),\;1\leq p<\infty$
denotes the space of measurable and integrable in power $p$
functions $f:\OO\rightarrow\RR$ with the norm
$$\|f\|_{L_p(\OO)}=\left(\int\limits_{\OO}|f(\vx)|^p\;d\vx\right)^{\frac 1
p},\;1\leq p<\infty.
$$

In order to introduce the set of functions that we will use as an
approximation tool, let us denote by
 $H(\OO)$ the set of harmonic functions in bounded domain $\OO\subset\RR^2$, i.e.
 $$H(\OO)=\left\{u(\vx)\in C^2(\OO):\;\Delta u(\vx)=0\right\},
$$
where $\Delta=\frac{\partial^2}{\partial
x_1^2}+\frac{\partial^2}{\partial x_2^2}$ is the Laplace operator.

Now we shall consider the set $S(\square_N)$ of continuous on
$\CI^2$ functions such that when restricted to any $\OO^N$
from the partition $\square_N$ are harmonic functions, i.e.
$$S(\square_N)=\left\{g\in C(\CI^2):\;g|_{\OO^N}=u_i|_{\OO^N},\;
u_i\in H(\OO^N)\;\forall\,\OO^N\in\square_N\right\}.
$$
From now on we shall
refer to the functions from set
 $S(\square_N)$ as {\it harmonic splines}.

Let now a function $f\in C^2(\CI^2)$ is given. We will approximate
it by a harmonic spline $s(f,\square_N)$ in such a way that their
values coincide on the boundary of elements $\OO^N$. In other words
we shall require the spline $s(f,\square_N)$ to ``interlinate'' the
given function on the boundary of all partition elements. Hence, the
function $s(f,\square_N)$ in each $\OO^N$ shall satisfy the equation
$$\Delta s(f,\square_N)=0
$$
and the boundary condition
$$s(f,\square_N)|_{\PO^N}=f|_{\PO^N}.
$$
In other words, function $s(f,\square_N)$ in each domain $\square_N$
must be a solution of Dirichlet problem for Laplace equation, which
is unique (see, for example, \cite{Riht})  and therefore harmonic
spline $s(f,\square_N)$ for each $f\in C^2(\CI^2)$ is well defined
on $\CI^2$.

The main goal of this paper is to study the minimal $L_p$-error ($1\leq p<\infty$) of approximation of a
given function $f$ by harmonic splines $s(f,\square_N)$ for all
possible partitions $\square_N$ which we shall denote by
$$R_N(f,L_p)=\inf_{\square_N}\|f-s(f,\square_N)\|_{L_p(\II^2)}.
$$
 The value of $R_N(f,L_p)$ we shall call the {\it optimal
error}. The partition $\widetilde{\square}_N$ on which the optimal
error is achieved will be called {\it optimal partition}, i.e.
$$R_N(f,L_p)=\|f-s(f,\widetilde{\square}_N)\|_{L_p(\II^2)}.
$$
For an arbitrary function
 $f$ it is impossible (except for some trivial cases) to find the optimal
partition $\widetilde{\square}_N$ and explicitly compute the value
of the optimal error. Therefore, we are interested in the following
two natural questions: first of all, how the optimal error behaves
as number of partition elements increases ($N\to \infty$), and
secondly, how to construct a
sequence of partitions $\{\square^*_N\}_{N=1}^{\infty}$, that is
{\it asymptotically optimal}, i.e.
$$\lim_{N\rightarrow\infty}\frac{\|f-s(f,\square^*_N)\|_{L_p(\II^2)}}{R_N(f,L_p)}=1.
$$

In order to state the main result of this paper, we need the
following notation. By $G_{\Omega}(\vx;\vv)$ with $\vx,\vv\in\Omega$
we shall denote the Green's function of the Dirichlet problem for
the domain $\Omega$ (for the detailed definition, see, for example, \cite{Riht}). In addition, denote by
\begin{equation}
I(\vx):=\int\limits_{\II^2}G_{\II^2}(\vx;\vv)d\vv,\;\vx\in\II^2.
\label{I}
\end{equation}

The main result of this paper is the following theorem.
\begin{theorem}\label{teorema1}
For an arbitrary function $f\in C^2(\CI^2)$ there exists a sequence
$\{\square^*_N\}_{N=1}^{\infty}$ and corresponding sequence of
splines $\{s(f,\square^*_N)\}_{N=1}^{\infty}$, such that
$$\lim_{N\to\infty}N\|f-s(f,\square^*_N)\|_{L_p(\II^2)}=\|I\|_{L_p(\II^2)}\|\Delta
f\|_{L_{\frac p {p+1}}(\II^2)}.
$$
\end{theorem}

{\bf Remark 1.} Note that the majority of the partition
elements from $\square^*_N$ are squares. This is rather restrictive
condition on the mesh. However, it is somewhat natural to expect
since we are using more information about the function than in the case of interpolation by polynomial splines.

{\bf Remark 2.} In the case when the domain of the function is
rectangle, the mesh elements must have the same shape as the domain,
i.e. the same width to length ratio.

Comparing this result with analogous results in the case of
approximating the given function from $C^2(\CI^2)$ by bilinear
polynomial splines (\cite{PhD,nasha}), we see that the order of the
error of approximation is the same. Therefore, depending on the
particular problem one could choose either polynomial or harmonic
splines, the latter providing interpolation on a larger set, but
allowing less flexible partitions.

\section{Idea of the proof and auxiliary results} \label{S3}

The main idea of the construction of each partition in the sequence
is to use first ``intermediate approximation'' of the given function
$f$ which needs to be approximated by splines. As an intermediate
approximation we shall use the piecewise function ``glued'' from
second degree Taylor polynomials of $f$. Then the intermediate
approximation is used to build the asymptotically optimal sequence of partitions of the
domain and, consequently, the sequence of corresponding harmonic splines.

The process of constructing this sequence of partitions consists of two steps:

\begin{enumerate}

\item
The domain $\CI^2$ is first divided into some number (small
comparing to $N$) of equal sizesubdomains, and on each instead of the original function
$f$ we consider its second degree Taylor polynomial constructed at,
say, center of the corresponding subdomain. The number of these
subdomains is chosen so that the absolute value of the difference
between the function and its Taylor polynomial is small enough on
each corresponding subdomain.

\item
Next, each subdomain is refined further, depending on the Taylor
polynomial (and therefore the original function) on this subdomain.
The final partition $\square^*_N$ will consist of the squares and
possibly some rectangles (their number is small comparing to $N$).
The total number of partition elements is determined from the
following two conditions: first of all, the total number of all the
elements in $\square^*_N$ is approximately equal to $N$, and
secondly, the global error of approximation of piecewise quadratic
function (consisting of Taylor polynomials built on the previous
step), by corresponding harmonic splines is minimal.

The problem of finding this error of approximation on ``smallest'
subdomain elements is reduced to approximating the function
$Q(\vx)=Ax_1^2+Bx_2^2$, where $A$ and $B$ are constants, since the
rest of the terms in second degree Taylor polynomial constitute a
harmonic function.

\end{enumerate}

The constructed in such a way partition $\square^*_N$ for each fixed
$N$ is used in the proof of the main Theorem \ref{teorema1}.

The idea of this construction based on intermediate approximation by
piecewise quadratic functions has been already used in papers
~\cite{us}, ~\cite{BBS}, ~\cite{PhD}, ~\cite{JAT}, ~\cite{nasha} by authors and co-authors to build
asymptotically optimal sequences of polynomial splines on
triangulations, rectangular partitions, and their generalizations in
various settings. More
on the history of adaptive approximation and asymptotically optimal
sequences of splines in the case of polynomial splines can be found
in ~\cite{BBS}, ~\cite{PhD}. 

Let us turn now to the auxiliary results that we shall need to prove the main theorem of this paper.

For a function $f\in C^2(\CI^2)$ we shall define a modulus of
continuity as follows
$$ \omega(f,\delta)=\sup
\{|f(\vx)-f(\vx')|:\;\; |x_1-x'_1|\leq \delta, |x_2-x'_2|\leq
\delta,\;\;\vx,\vx'\in \CI^2\}.
$$
Then we consider
\begin{equation}
\omega(\delta)=\max\left\{\omega(f_{xx},\delta),\omega(f_{yy},\delta),\omega(f_{xy},\delta)\right\}.\label{omega}
\end{equation}
\begin{lemma}\label{lemma1}
Let function $f\in C^2(\CI^2)$ be given and let
$P_2(\vx)=P_2(f;\vx;\vx_0)$ be its second order Taylor polynomial
taken at the point $\vx_0$, that is the center of some square
$D_h\subset \CI^2$ with side length $h$. Then
\begin{equation}
|f(\vx)-P_2(\vx)|\leq \frac{h^2}{2}\omega \left(\frac{h}{2} \right
), \qquad \vx\in D_h.
\end{equation}
\end{lemma}

For the proof of this lemma see, for example, ~\cite{us}.

To state the next lemma we shall need the following notation.

Let $G_{\Omega}(\vx;\vv)$ be a Green's function of the interior
Dirichlet problem for bounded region $\Omega\subset\RR^2$. By
$\Omega+\vd$ we denote the region obtained by the shift of the
region $\Omega$ by a vector $\vd\in\RR^2$, and by $\alpha\Omega$ we
denote the region obtained by scaling $\Omega$ by coefficient
$\alpha$.
% -- область, которая является
%гомотетией области $\Omega$ с центром в точке $a=(a_1,a_2)$ и
%коэффициентом $\alpha$.

\begin{lemma}\label{lemma2} The following properties of the Green's function hold true:
\begin{equation}
G_{\alpha\Omega}(\vx;\vv)=
%\frac 1
%{\alpha^{n-2}}
G_{\Omega}\left(\frac \vx {\alpha};\frac \vv
{\alpha}\right),\;\vx\in\alpha\Omega,\vv\in\overline{\alpha\Omega},\label{gr1}
\end{equation}
\begin{equation}
G_{\Omega+d}(\vx;\vv)=G_{\Omega}(\vx-\vd;\vv-\vd),\;\vx\in\Omega+\vd,\vv\in\overline{\Omega+\vd}.\label{gr2}
\end{equation}
\end{lemma}

{\bf Proof.} In the proof of this lemma we shall use the well-known
(see, for instance, ~\cite{Riht}) formula for the Green's
function of the interior Dirichlet problem for the bounded region
$\Omega$. For
 $\vx\in\Omega,\;\vv\in\CO$ the Green's function can be written as
 $$G_{\Omega}(\vx;\vv)=\frac 1 {2\pi}\ln\frac 1
{|\vx-\vv|}+g(\vx;\vv),
$$
where function $g(\vx;\vv)$ is harmonic in both arguments on $\OO$,
continuous in $\vv$ on $\CO$, and is chosen so that the Green's
function has zero value on the boundary.

For $\vx\in\alpha\Omega,\;\vv\in\overline{\alpha\Omega}$
$$G_{\Omega}\left(\frac {\vx} {\alpha};\frac {\vv} {\alpha}\right)=\frac 1
{2\pi}\ln\frac 1 {\left|\frac \vx {\alpha}-\frac \vv
{\alpha}\right|}+g\left(\frac \vx {\alpha};\frac \vv
{\alpha}\right)=\frac 1 {2\pi}\ln\frac {\alpha}
{\left|\vx-\vv\right|}+g\left(\frac \vx {\alpha};\frac \vv
{\alpha}\right)=
$$
$$=\frac 1 {2\pi}\left(\ln\alpha+\ln\frac 1
{\left|\vx-\vv\right|}\right)+g\left(\frac \vx {\alpha};\frac \vv
{\alpha}\right)=\frac 1 {2\pi}\ln\frac 1
{\left|\vx-\vv\right|}+\tilde g(\vx;\vv),
$$
where function $\tilde g(\vx;\vv)=g\left(\frac \vx {\alpha};\frac
\vv {\alpha}\right)+\frac 1 {2\pi}\ln\alpha$ is clearly harmonic in
both arguments in $\alpha\Omega$ and continuous in $\vv$ on
$\overline{\alpha\Omega}$. Besides that, since the boundary of the
region $\Omega$ is mapped to the boundary of the region
$\alpha\Omega$, and the Green's function has zero value on
$\partial\Omega$, then the expression above will be zero on
$\partial(\alpha\Omega)$. Hence, (\ref{gr1}) is proved.

Let us consider next the shift of the region $\Omega$. For
$\vx\in{\Omega+\vd},\;\vv\in\overline{\Omega+\vd}$
$$G_{\Omega}(\vx-\vd;\vv-\vd)=\frac 1 {2\pi}\ln\frac 1 {|(\vx-\vd)-(\vv-\vd)|}+g(\vx-\vd;\vv-\vd)=
$$
$$=\frac 1 {2\pi}\ln\frac 1 {|\vx-\vv|}+\tilde g(\vx;\vv).
$$
Since all the necessary conditions for the function
 $\tilde
g(\vx;\vv)=g(\vx-\vd;\vv-\vd)$ are satisfied, then the obtained expression is the Green's function for the region $\Omega+\vd$, and hence (\ref{gr2}) is proved. $\square$

Next we shall consider a square region $\Omega\subset\RR^2$ with sides parallel to the coordinate axis, and on which the following functions is defined
\begin{equation}
Q(\vx)=Ax_1^2+Bx_2^2,\label{Q}
\end{equation}
where $A$ and $B$ are constants.

We shall find the $L_p$-error of approximation of the function
 $Q(\vx)$ by a harmonic function $u(Q;\vx)$ on $\Omega$ such that the values of $u(Q;\vx)$ coincide with the function  $Q(\vx)$ on the boundary of $\Omega$.
 In other words, for the function $u(Q;\vx)$ the following conditions must be satisfied
\begin{equation}
\Delta u(Q;\vx)=0,\;\vx\in\Omega\label{Usl_1}
\end{equation}
and
\begin{equation}
u(Q;\vx)|_{\partial\Omega}=Q(\vx)|_{\partial\Omega}.\label{Usl_2}
\end{equation}

\begin{lemma}\label{lemma3} For the quadratic function $Q(\vx)=Ax_1^2+Bx_2^2$
such that $\Delta Q(\vx)=2(A+B)\neq0$, we have\begin{equation}
\|Q-u(Q)\|_{L_p(\Omega)}=2|A+B|\cdot|\Omega|^{1+\frac 1
p}\|I\|_{L_p(\II^2)},
\end{equation}
where $|\Omega|$ denotes the area of $\Omega$, and $I$
is defined in (\ref{I}).
\end{lemma}

{\bf Proof.} We shall use the following representation of an arbitrary function $u(\vx)\in C^2(\Omega)\bigcap C(\overline{\Omega})$ on the region $\Omega\subset\RR^2$ with smooth enough boundary:
\begin{equation}
u(\vx)=-\int\limits_{\partial\Omega}\frac{\partial
G_{\Omega}(\vx;\vv)}{\partial\bar n}u(\vv)d\vv
-\int\limits_{\Omega}G_{\Omega}(\vx;\vv)\Delta
u(\vv)d\vv,\;\vx\in\Omega,\label{osn}
\end{equation}
where $\frac{\partial G_{\Omega}(\vx;\vv)}{\partial\bar n}$ is its
partial derivative in the outer normal direction to $\partial\Omega$
(see, for example, ~\cite{Riht}).

Then for the function $Q(\vx)$, taking into account (\ref{osn}) and
the fact that $\Delta Q(\vx)=2(A+B)$, we have
\begin{equation}
Q(\vx)=-\int\limits_{\partial\Omega}\frac{\partial
G_{\Omega}(\vx;\vv)}{\partial\bar n}Q(\vv)d\vv
-2(A+B)\int\limits_{\Omega}G_{\Omega}(\vx;\vv)d\vv,\;\vx\in\Omega.\label{Q_osn}
\end{equation}

In addition, with the help of (\ref{osn}), the solution of the
problem (\ref{Usl_1}) -- (\ref{Usl_2}) can be written as
\begin{equation}
u(Q;\vx)=-\int\limits_{\partial\Omega}\frac{\partial
G_{\Omega}(\vx;\vv)}{\partial\bar
n}Q(\vv)d\vv,\;\vx\in\Omega.\label{u_osn}
\end{equation}

From (\ref{Q_osn}), (\ref{u_osn}), and the fact that the Green's function is nonnegative, it follows that
$$
|Q(\vx)-u(Q;\vx)|=2|A+B|\int\limits_{\Omega}G_{\Omega}(\vx;\vv)d\vv,\;\vx\in\Omega,
$$
and, hence,
\begin{equation}
\|Q-u(Q)\|_{L_p(\Omega)}=2|A+B|\left(\int\limits_{\Omega}\left(\int\limits_{\Omega}G_{\Omega}(\vx;\vv)d\vv\right)^pd\vx\right)^{\frac
1 p}.\label{L_p}
\end{equation}
Our next goal is to transform the integral in the right-hand side of
(\ref{L_p}) so that it explicitly depends on the area of $\Omega$.
In order to do so we need the result of Lemma~\ref{lemma2}.

Let $\widetilde{\II}^2$ be the square unit region such that when centered at one of its vertices and stretched by $\alpha$ it becomes
$\Omega$. Then the volumes of these regions are related by $|\Omega|=\alpha^2|\widetilde{\II}^2|$, from where we have
$\alpha=\sqrt{|\Omega|}$.

Therefore, we have the following equality
\begin{equation}
\left(\int\limits_{\Omega}\left(\int\limits_{\Omega}G_{\Omega}(\vx;\vv)d\vv\right)^pd\vx\right)^{\frac
1 p}
=\left(\int\limits_{\alpha\widetilde{\II}^2}\left(\int\limits_{\alpha\widetilde{\II}^2}G_{\widetilde{\II}^2}\left(\frac
\vx {\alpha};\frac \vv
{\alpha}\right)d\vv\right)^pd\vx\right)^{\frac 1 p}.\nonumber
\end{equation}

Changing variables we arrive at
\begin{equation}
\left(\int\limits_{\Omega}\left(\int\limits_{\Omega}G_{\Omega}(\vx;\vv)d\vv\right)^pd\vx\right)^{\frac
1 p} =\alpha^{2+\frac 2
p}\left(\int\limits_{\widetilde{\II}^2}\left(\int\limits_{\widetilde{\II}^2}G_{\widetilde{\II}^2}(\vx;\vv)d\vv\right)^pd\vx\right)^{\frac
1 p}.\label{int}
\end{equation}

It is left to notice that due to property (\ref{gr2}) of the Green's function, the last equality will not change after the shift of the region
$\widetilde{\II}^2$ by an arbitrary vector from $\RR^2$. Therefore, as
 $\widetilde{\II}^2$ we may take an arbitrary rectangular region of unit area with sides parallel to the coordinate axis, in particular, we may take $\II^2$.

Taking into account (\ref{int}), notation (\ref{I}), and having recalled that for the given region $\Omega$ the value
$\alpha=\sqrt{|\Omega|}$, equality (\ref{L_p}) becomes
$$
\|Q-u(Q)\|_{L_p(\Omega)}=2|A+B|\cdot|\Omega|^{1+\frac 1
p}\|I\|_{L_p(\II^2)}.
$$
$\square$

Now let us turn to the proof of Theorem~\ref{teorema1}.

\section {Proof of the main result} \label{S5}

We shall begin the proof by constructing (for a given function
 $f\in
C^2(\CI^2)$ and fixed number $N$) the partition $\square^*_N$ of
$\CI^2$ such that the sequence $\{\square^*_N\}_{N=1}^{\infty}$ is asymptotically optimal. The idea of this
construction is given in Section \ref{S3}.

First of all, we determine the number of the elements of the intermediate partition. To this end, for an arbitrary fixed
$\varepsilon\in (0,1)$ we set
$$
m_N:=\min\left \{ m\in \NN : \;\; \frac1
2\left(\frac1{m}\right)^2\omega \left(\frac{1}{2m}\right)\leq
\frac{\varepsilon}{N}\right\}, $$
where $\omega(\delta)$ was defined
in (\ref{omega}). It is clear that $m_N\rightarrow\infty$ as
$N\rightarrow\infty$.

We shall take $m_N^2$ as the number of elements of intermediate
partition, and we shall subdivide the square $\CI^2$ into equal
squares $D_N^l,\;\;l=1,...,m_N^2$ with the side length
$\frac{1}{m_N}$. Let us show that the number of elements in
intermediate partition is indeed small comparing to the total number
of elements $N$, i.e.
\begin{equation}
m_N^2=o(N),\;\;N\rightarrow\infty.\label{moN}
\end{equation}
In order to do so let us consider $\frac{N}{m_N^2}$:
\begin{multline}
\frac{N}{m_N^2}=\frac N 2
\frac1{(m_N-1)^2}\omega\left(\frac1{2(m_N-1)}\right)\frac{2(m_N-1)^2}{m^2_N\omega\left(\frac{1}{2(m_N-1)}\right)}\\
\geq\varepsilon\frac{2(m_N-1)^2}{m^2_N\omega\left(\frac{1}{2(m_N-1)}\right)}.
\end{multline}
Taking into account $\left(\frac{m_N-1}{m_N}\right)\rightarrow1$ and
$\omega\left(\frac1{2m_N}\right)\rightarrow0$ as
$N\rightarrow\infty$ we obtain $\frac{N}{m_N^2}\rightarrow\infty$ as
$N\rightarrow\infty$, and hence (\ref{moN}) holds true.

In order to obtain the intermediate estimates we define the
functions $f_N(\vx)$ and $Q_N(\vx)$, using the notation
$P_2(f;\vx;\vx_l)$ for Taylor polynomial of second degree for
$f(\vx)$ at the center $\vx_l$ of each square $D_l^N$,
as follows:
\begin{enumerate}
\item For $\vx\in D_1^N$ set
$$f_N(\vx)=P_2(f;\vx;\vx_1)$$
and
$$Q_N(\vx)=\frac 1 2 \frac{\partial^2f}{\partial x_1^2}(\vx_1)x_1^2+\frac 1 2 \frac{\partial^2f}{\partial x_2^2}(\vx_1)x_2^2.
$$
\item When $1<l\leq m^2_N$ for $\vx\in D_l^N\setminus\cup_{i=1}^{l-1}\partial D_i^N$ we set
$$f_N(\vx)=P_2(f;\vx;\vx_l)$$
and
$$Q_N(\vx)=\frac 1 2 \frac{\partial^2f}{\partial x_1^2}(\vx_l)x_1^2
+\frac 1 2 \frac{\partial^2f}{\partial x_2^2}(\vx_l)x_2^2.
$$
\end{enumerate}

Next we shall find the number of elements to additionally subdivide
each $D_N^l,\;\;l=1,...,m_N^2$. For that we shall cover the square
$D_N^l$ by a ``mesh'' consisting of squares of {the same} fixed
area, fixing, for convenience, one of the mesh vertices with one of
the vertices of the square $D_N^l$. The intersection of $D_N^l$ with
this mesh will provide the necessary subdivision of the square. It
will consist of the squares of the original mesh as well as,
possibly, some rectangles along the boundary of $D_N^l$.

We shall consider two cases: when the partition consists only of squares and when the partition contains some rectangles along the boundary.

First, let us assume that partition of $D_N^l$ consists of squares only. We shall find the area of the squares next.
%In order to find the area of the squares in the mesh for the given function
% $f(\vx)$ and fixed $N$, предположим, что наше
%разбиение состоит только из квадратов.
The number of the squares in the partition of $D_N^l$ we shall denote by $\tilde{n}_l^2$. The elements itself we shall denote by $\OO_N^{l,i}$,
$i=1,...,\tilde n_l^2,\;l=1,...,m_N^2$. Then
$\bigcup\limits_{l=1}^{m_N^2}\bigcup\limits_{i=1}^{\tilde n_l^2}\OO_N^{l,i}$
gives the partition of $\CI^2$ (for fixed
$N$), and hence $\sum\limits_{l=1}^{m_N^2}\tilde n_l^2=N$.

To determine $\tilde n_l^2$ for each $D_N^l$ we shall minimize the
value of the error of approximation of $f_N(\vx)$ by the
corresponding harmonic spline $s(f_N,\square^*_N)$ on $\CI^2$. The
needed error can be written in terms of the errors on each element
of the partition as
\begin{equation}
\|f_N-s(f_N,\square^*_N)\|^p_{L_p(\II^2)}=\sum_{l=1}^{m_N^2}\sum_{i=1}^{n_l^2}\|f_N-s(f_N,\square^*_N)\|^p_{L_p(\OO_N^{l,i})}.\label{pogrNaI}
\end{equation}
Therefore, we need to find the error of approximation of $f_N(\vx)$
on each element
 $\OO_N^{l,i}$.
Since on $\OO_N^{l,i}$ the difference $f_N(\vx)-Q_N(\vx)$ is a
harmonic function, using Lemma~\ref{lemma3}, we have
\begin{multline}
\|f_N-s(f_N,\square^*_N)\|^p_{L_p(\OO_N^{l,i})}=\|Q_N-s(f_N,\square^*_N)\|^p_{L_p(\OO_N^{l,i})}\\
=\left|\frac{\partial^2 f}{\partial x^2}(\vx_l)+\frac{\partial^2
f}{\partial y^2}(\vx_l)\right|^p|\OO_N^{l,i}|^{p+1}
\|I\|^p_{L_p(\II^2)}=M^p_{f,N}(\vx_l)\frac{1}{(m_N\tilde
n_l)^{2(p+1)}}\|I\|^p_{L_p(\II^2)},\label{pogrNaOm}
\end{multline}
where we used the notation
\begin{equation}
M_{f,N}(\vx_l)=\left|\frac{\partial^2 f}{\partial
x^2}(\vx_l)+\frac{\partial^2 f}{\partial
y^2}(\vx_l)\right|.\label{defM}
\end{equation}
Taking into account
(\ref{pogrNaOm}), the equality (\ref{pogrNaI}) becomes
$$
\|f_N-s(f_N,\square^*_N)\|^p_{L_p(\II^2)}=\sum_{l=1}^{m_N^2}M^p_{f,N}(\vx_l)\frac{1}{(m_N\tilde
n_l)^{2(p+1)}}\|I\|^p_{L_p(\II^2)}\tilde n_l^2,
$$
from where we finally arrive at
\begin{equation}
\|f_N-s(f_N,\square^*_N)\|^p_{L_p(\II^2)}=\sum_{l=1}^{m_N^2}M^p_{f,N}(\vx_l)\frac{1}{(m_N)^{2(p+1)}(\tilde
n_l)^{2p}}\|I\|^p_{L_p(\II^2)}.\label{pogrNaIFin}
\end{equation}
Using the method of Lagrange multipliers, we shall minimize the right-hand side of
 (\ref{pogrNaIFin}) under constraint
\begin{equation}
\sum\limits_{l=1}^{m_N^2}\tilde n_l^2=N.\label{usl_sum_n}
\end{equation}
For that let us consider the function
\begin{equation}
\Phi(\tilde n_1,..,\tilde
n_{m_N^2},\lambda)=\sum_{l=1}^{m_N^2}M^p_{f,N}(\vx_l)\frac{1}{(m_N)^{2(p+1)}(\tilde
n_l)^{2p}}\|I\|^p_{L_p(\II^2)}
+\lambda\sum\limits_{l=1}^{m_N^2}\tilde n_l^2,\label{funcPhi}
\end{equation}
where $\lambda$ is the Lagrange multiplier.

To find the critical points of the function $\Phi(\tilde
n_1,..,\tilde n_{m_N^2},\lambda)$ we shall consider the system
\begin{equation}\left\{
\begin{array}{lll}
\frac{\partial}{\partial \tilde
n_l}\left(\sum\limits_{l=1}^{m_N^2}M^p_{f,N}(\vx_l)\frac{1}{(m_N)^{2(p+1)}(\tilde
n_l)^{2p}}\|I\|^p_{L_p(\II^2)}
+\lambda\sum\limits_{l=1}^{m_N^2}\tilde
n_l^2\right)=0,\;l=1,...,m_N^2,\cr \sum\limits_{l=1}^{m_N^2}\tilde
n_l^2=N.\cr
\end{array}\right.\label{sistema}
\end{equation}
Solving it we obtain $\tilde n_l^2$:
\begin{equation}
\tilde
n_l^2=\frac{NM_{f,N}^{\frac{p}{p+1}}(\vx_l)}{
\sum\limits_{i=1}^{m_N^2}M_{f,N}^{\frac{p}{p+1}}(\vx_i)},
\;l=1,...,m_N^2.\label{nMin}
\end{equation}

Since for the solutions of the system (\ref{sistema}) the value
$\textrm{d}^2\Phi>0$, then (\ref{nMin}) indeed provide the minimal
value of the total error of approximation of the function $f_N(\vx)$
by harmonic spline $s(f_N,\square^*_N)$ on the square $\CI^2$.

As the area of every square of the mesh
$\OO\subset D_N^l$ for the fixed $N$ we shall take
\begin{equation}
|\OO|=\frac{1}{(m_N\tilde{n}_l)^2},\;l=1,...,m_N^2.\label{setka}
\end{equation}
Later we shall need the following estimate for the values $\tilde{n}_l$
\begin{equation}
\tilde{n}_l\geq\frac{\sqrt{N}\min\limits_{\vx\in\II^2}\{M_{f,N}^{\frac{p}{2(p+1)}}(\vx)\}}{m_N\|\Delta
f\|_{L_{\infty}(\II^2)}^{\frac{p}{2(p+1)}}}=C\frac{\sqrt{N}}{m_N}>0,\label{ocenka_n}
\end{equation}
where $C$ is independent of $N$ constant.

%\tc{Zametim, chto dlya togo, chtoby okonchatel'noe razbienie oblasti
%opredeleniya sostoyalo tol'ko iz kvadratov, velichiny, opredelyaemye
%ravenstvom \ref{nMin} dlya vsex $l=1,...,m_N^2$ dolzhny yavlyat'sya
%kvadratami natural'nyx chisel. V etom sluchae my poluchim
%optimal'noe razbienie.}

For the final partition of the domain to consist of only squares, the values defined by (\ref{nMin}) for all $l=1,...,m_N^2$ must be squares of natural numbers. In this case we will have optimal partition. Otherwise, the final partition will contain both squares and rectangles, and will be optimal only asymptotically.

From now on we shall assume that the partition of
 $D_N^l$ consists of both squares of the {mesh} and some rectangles along the boundary. By
$R_N^l$ we shall denote the set of all squares, and by  $\widetilde R_N^l$ - the set of rectangles from the partition of $D_N^l$. In addition, let $n_l^2=n_l^2(N)$ denote the number of squares from
$R_N^l$. Then the number of rectangles from $\widetilde R_N^l$
is clearly equal to $2n_l+1$.

It is also clear that for all $l=1,...,m_N^2$ we have
\begin{equation}
(\tilde{n}_l-1)^2<n_l^2\leq\tilde{n}_l^2.\label{def_n}
\end{equation}
The number of elements of the partition of
$D_N^l$ (in the case when rectangles are present) is $(n_l+1)^2$. Let us show that the total number of the partition elements $\sum\limits_{l=1}^{m_N^2}(n_l+1)^2$
and $N$ are values of the same order as $N\rightarrow\infty$.
To that end, let us estimate the values $\sum\limits_{l=1}^{m_N^2}n_l^2$ and $\sum\limits_{l=1}^{m_N^2}n_l$.
Using inequality (\ref{def_n}) and condition (\ref{usl_sum_n})
we have
\begin{equation}
\sum\limits_{l=1}^{m_N^2}n_l^2\leq\sum\limits_{l=1}^{m_N^2}\tilde{n}_l^2=N,\label{ocenka1}
\end{equation}
\begin{equation}
\sum\limits_{l=1}^{m_N^2}n_l^2\geq\sum\limits_{l=1}^{m_N^2}\tilde{n}_l^2-2
\sum\limits_{l=1}^{m_N^2}\tilde{n}_l=N-2\sum\limits_{l=1}^{m_N^2}\tilde{n}_l.
\end{equation}
Using Holder inequality, we obtain
\begin{equation}
\sum\limits_{l=1}^{m_N^2}n_l\leq\sum\limits_{l=1}^{m_N^2}\tilde{n}_l
\leq
m_N\left(\sum\limits_{l=1}^{m_N^2}\tilde{n}_l^2\right)^{\frac{1}{2}}=m_N\sqrt
N.\label{ocenka3}
\end{equation}
Taking into account estimates (\ref{ocenka1}) -- (\ref{ocenka3}), we have the following double inequality
$$
N-2m_N\sqrt N\leq\sum\limits_{l=1}^{m_N^2}(n_l+1)^2\leq N+2m_N\sqrt
N+m_N^2.$$
Dividing both sides by
 $N$, taking the limit as $N\rightarrow\infty$ and using (\ref{moN}), we obtain
$$\lim\limits_{N\rightarrow\infty}\frac{\sum\limits_{l=1}^{m_N^2}(n_l+1)^2}{N}=1.
$$
Therefore, the number of elements of the constructed partition has order $N$ as $N\to \infty$.

Next we shall verify that the constructed in such a way sequence of partitions will be asymptotically optimal.

%Начнем с получения оценки сверху для величины погрешности
%приближения заданной функции $f(x)$ гармоническим сплайном
%$s(f,\square_N)$. Для этого нам понадобится ранее определенная
%функция $f_N(x)$ и приближающий ее гармонический сплайн
%$s(f_N,\square_N)$. Имеет место следующее неравенство

Using the triangle inequality we have
\begin{multline}
\|f-s(f,\square_N)\|_{L_p(\II^2)}\leq\|f-f_N\|_{L_p(\II^2)}+\|s(f,\square_N)-s(f_N,\square_N)\|_{L_p(\II^2)}\\
+\|f_N-s(f_N,\square_N)\|_{L_p(\II^2)}.\label{ocenka_sverxu}
\end{multline}
We shall now estimate each term on the right.

For the first and second terms we shall first obtain the estimates on an arbitrary square
 $D_N^l$ of intermediate partition in terms of $\|f-f_{N}\|_{L_{\infty}(D_N^l)}$.
Since $|D_N^l|=\frac{1}{m_N^2}$, then the following is true
\begin{equation}
\|f-f_{N}\|^p_{L_p(D_N^l)}\leq\frac{1}{m_N^2}\|f-f_{N}\|^p_{L_{\infty}(D_N^l)}.\label{normRaznF}
\end{equation}
To estimate the difference
$\|s(f,\square^*_N)-s(f_{N},\square^*_N)\|_{L_p(D_N^l)}$ we shall consider an arbitrary element of the partition $\OO\subset D_N^l$. According to the definition of a harmonic spline, the difference of functions $s(f,\square^*_N)(\vx)-s(f_{N},\square^*_N)(\vx)$ is a harmonic function on $\OO$ and its values on $\partial\OO$ coincide with the values of the difference $f(\vx)-f_{N}(\vx)$. Therefore, using representation (\ref{osn}), we obtain
$$
s(f,\square^*_N)(\vx)-s(f_{N},\square^*_N)(\vx)=-\int\limits_{\partial\OO}\frac{\partial
G_{\Omega}(\vx;\vv)}{\partial\bar
n}(f(\vv)-f_{N}(\vv))d\vv,\;\vx\in\OO.
$$
Then
$$|s(f,\square^*_N)(\vx)-s(f_{N},\square^*_N)(\vx)|\leq\|f-f_{N}\|_{L_{\infty}(\OO)}\int\limits_{\partial\OO}\frac{\partial
G_{\Omega}(\vx;\vv)}{\partial\bar n}d\vv,\;\vx\in\OO.
$$
Taking into account the uniqueness of the solution (\ref{osn}) (see,
for instance, ~\cite{Riht}), we have
$$
\int\limits_{\partial\OO}\frac{\partial
G_{\Omega}(\vx;\vv)}{\partial\bar n}d\vv=1,
$$
and, therefore, for all $\vx\in D_N^l$
$$|s(f,\square^*_N)(\vx)-s(f_{N},\square^*_N)(\vx)|\leq\|f-f_{N}\|_{L_{\infty}(D_N^l)},
$$
which implies
\begin{equation}
\|s(f,\square^*_N)-s(f_{N},\square^*_N)\|^p_{L_p(D_N^l)}\leq\frac{1}{m_N^2}\|f-f_{N}\|^p_{L_{\infty}(D_N^l)}.\label{normRaznS}
\end{equation}
Using now Lemma~\ref{lemma1} and definition of $m_N$, we obtain
$$\|f-f_{N}\|_{L_{\infty}(D_l^N)}\leq\frac{1}{2m_N^2}\omega\left(\frac{1}{2m_N}\right)\leq\frac{\varepsilon}{N}.
$$
The last inequality, together with
 (\ref{normRaznS}) and
(\ref{normRaznF}), imply the needed estimates:
\begin{equation}
\|f-f_{N}\|_{L_p(\II^2)}\leq\frac{\varepsilon}{N},\label{normRaznFFin}
\end{equation}
\begin{equation}
\|s(f,\square^*_N)-s(f_{N},\square^*_N)\|_{L_p(\II^2)}\leq\frac{\varepsilon}{N}.\label{normRaznSFin}
\end{equation}

Next let us consider the third term in (\ref{ocenka_sverxu})
\begin{multline}
\|f_N-s(f_N,\square_N)\|^p_{L_p(\II^2)}=\sum\limits_{l=1}^{m_N^2}\|f_N-s(f_N,\square_N)\|^p_{L_p(D_N^l)}\\
=\sum\limits_{l=1}^{m_N^2}\left(\sum\limits_{\OO\in
R_N^l}\|f_N-s(f_N,\square_N)\|^p_{L_p(\OO)}+\sum\limits_{\widetilde\OO\in
\widetilde
R_N^l}\|f_N-s(f_N,\square_N)\|^p_{L_p(\widetilde\OO)}\right).\label{3slag}
\end{multline}
Since for each rectangle $\widetilde{\OO}\in \widetilde
R_N^l$ there exists such an element of the mesh $\OO$ that
$\widetilde{\OO}\subset\OO$, then
$$\|f_N-s(f_N,\square_N)\|^p_{L_p(\widetilde\OO)}\leq\|f_N-s(f_N,\square_N)\|^p_{L_p(\OO)}.
$$
Therefore, using equality (\ref{pogrNaOm}) for expression
(\ref{3slag}) we obtain
$$\|f_N-s(f_N,\square_N)\|^p_{L_p(\II^2)}\leq\sum\limits_{l=1}^{m_N^2}M^p_{f,N}(\vx_l)\frac{1}{(m_N\tilde
n_l)^{2(p+1)}}\|I\|^p_{L_p(\II^2)}(n_l+1)^2.
$$
Taking into account the estimate (\ref{def_n}) for the number of elements $n_l$, we have
\begin{multline}
\|f_N-s(f_N,\square_N)\|^p_{L_p(\II^2)}\leq\frac{\|I\|^p_{L_p(\II^2)}}{m_N^{2(p+1)}}
\sum\limits_{l=1}^{m_N^2}M^p_{f,N}(\vx_l)\frac{(\tilde
n_l+1)^2}{(\tilde n_l)^{2(p+1)}}\\
=\frac{\|I\|^p_{L_p(\II^2)}}{m_N^{2(p+1)}}
\sum\limits_{l=1}^{m_N^2}M^p_{f,N}(\vx_l)\frac{1}{\tilde{n}^{2p}_l}\left(1+\frac{2}{\tilde{n}_l}+\frac{1}{\tilde{n}^2_l}\right).\nonumber
\end{multline}

We shall estimate $\tilde{n}^{2p}_l$ using
(\ref{nMin}), and for $\tilde{n}_l$ and $\tilde{n}^2_l$
we shall use (\ref{ocenka_n})
\begin{multline}
\|f_N-s(f_N,\square_N)\|^p_{L_p(\II^2)}\\
\leq\frac{\|I\|^p_{L_p(\II^2)}}{m_N^{2(p+1)}}
\sum\limits_{l=1}^{m_N^2}M^p_{f,N}(\vx_l)\frac{\left(
\sum\limits_{i=1}^{m_N^2}M_{f,N}^{\frac{p}{p+1}}(\vx_i)\right)^p}{N^pM_{f,N}^{\frac{p^2}{p+1}}(\vx_l)}
\left(1+\frac{2m_N}{\sqrt N}+\frac{m_N^2}{N}\right)\\
=\frac{\|I\|^p_{L_p(\II^2)}}{N^p}\left(\frac{1}{m_N^2}
\sum\limits_{i=1}^{m_N^2}M_{f,N}^{\frac{p}{p+1}}(\vx_i)\right)^{p+1}\left(1+\frac{2m_N}{\sqrt
N}+\frac{m_N^2}{N}\right).
\end{multline}

Since
$$\frac{1}{m_N^2}
\sum\limits_{i=1}^{m_N^2}M_{f,N}^{\frac{p}{p+1}}(\vx_i)\rightarrow\int\limits_{\II^2}\Delta
f^{\frac{p}{p+1}}(\vx)d\vx,\;N\rightarrow\infty,
$$
then for large enough $N$
$$\|f_N-s(f_N,\square_N)\|^p_{L_p(\II^2)}
\leq\frac{\|I\|^p_{L_p(\II^2)}}{N^p}\left((1+\varepsilon)\int\limits_{\II^2}\Delta
f^{\frac{p}{p+1}}(\vx)d\vx\right)^{p+1}\left(1+\frac{2m_N}{\sqrt
N}+\frac{m_N^2}{N}\right).
$$
Going back to (\ref{ocenka_sverxu}), and taking into account
(\ref{normRaznFFin}) and (\ref{normRaznSFin}), we obtain
$$\|f-s(f_,\square_N)\|_{L_p(\II^2)}
\leq\frac{2\varepsilon}{N}+\frac{\|I\|_{L_p(\II^2)}}{N}\left((1+\varepsilon)\int\limits_{\II^2}\Delta
f^{\frac{p}{p+1}}(\vx)d\vx\right)^{\frac{p+1}p}\left(1+\frac{2m_N}{\sqrt
N}+\frac{m_N^2}{N}\right)^{\frac 1 p}.
$$
Hence, using
(\ref{moN}), we have
$$\limsup_{N\rightarrow\infty}N\|f-s(f,\square_N)\|_{L_p(\II^2)}=2\varepsilon+\|I\|_{L_p(\II^2)}\left((1+\varepsilon)\int\limits_{\II^2}\Delta
f^{\frac{p}{p+1}}(\vx)d\vx\right)^{\frac{p+1}p}.
$$
Finally, since $\varepsilon\in(0,1)$ is arbitrary, then switching to
the limit as
 $\varepsilon\rightarrow +0$,
we obtain
$$\limsup_{N\rightarrow\infty}N\|f-s(f,\square_N)\|_{L_p(\II^2)}=\|I\|_{L_p(\II^2)}\|\Delta
f\|_{L_{\frac{p}{p+1}}(\II^2)}.
$$
In order to obtain the estimate from below we shall use again the
function $f_N(\vx)$ for intermediate approximation. We shall use the
triangle inequality for
 $\|f_{N}-s(f_{N},\square^*_N)\|_{L_p(\II^2)}$ as follows
\begin{multline}
\|f-s(f,\square^*_N)\|_{L_p(\II^2)}\geq\|f_{N}-s(f_{N},\square^*_N)\|_{L_p(\II^2)}
-\|f-f_{N}\|_{L_p(\II^2)}\\
-\|s(f,\square^*_N)-s(f_{N},\square^*_N)\|_{L_p(\II^2)}.\label{lowEst}
\end{multline}
The estimate from below for the first term can be obtained very similarly to the way we obtained the estimate from above for this term. The result of it (for large enough $N$) will be
$$\|f_N-s(f_N,\square^*_N)\|_{L_p(\II^2)}>\frac{\|I\|_{L_p(\II^2)}}{N}\left((1-\varepsilon)\int\limits_{\II^2}|\Delta
f(\vx)|^{\frac{p}{p+1}}d\vx\right)^{\frac{p+1}{p}}\left(1-\frac{2m_N}{\sqrt
N}\right)^{\frac{1}{p}}.
$$
Now we may rewrite (\ref{lowEst}), using for second and third terms the inequalities
(\ref{normRaznFFin}) and (\ref{normRaznSFin})
$$\|f-s(f,\square^*_N)\|_{L_p(\II^2)}>\frac{\|I\|_{L_p(\II^2)}}{N}\left((1-\varepsilon)\int\limits_{\II^2}|\Delta
f(\vx)|^{\frac{p}{p+1}}d\vx\right)^{\frac{p+1}{p}}\left(1-\frac{2m_N}{\sqrt
N}\right)^{\frac{1}{p}}-\frac{2\varepsilon}{N}.
$$
Multiplying both sides by $N$, switching to the limit as
$N\rightarrow\infty$, and taking into account (\ref{moN})
$$\liminf_{N\rightarrow\infty}N\|f-s(f,\square^*_N)\|_{L_p(\II^2)}\geq\|I\|_{L_p(\II^2)}\left((1-\varepsilon)\int\limits_{\II^2}|\Delta
f(\vx)|^{\frac{p}{p+1}}d\vx\right)^{\frac{p+1}{p}}-{2\varepsilon}.
$$
Since $\varepsilon$ is arbitrary, we arrive at
$$\liminf_{N\rightarrow\infty}N\|f-s(f,\square^*_N)\|_{L_p(\II^2)}\geq\|I\|_{L_p(\II^2)}\|\Delta
f\|_{L_{\frac{p}{p+1}}(\II^2)}.
$$
Together with the estimate from above it completes the proof of
Theorem~\ref{teorema1}. $\square$

\begin{acknowledgements}
 Authors would like to thank Professor V. Babenko for his advice and guidance during work on this project.
\end{acknowledgements}

\end{document}